\documentclass[95pt]{article}

\usepackage{verbatim}
\usepackage{hyperref}
\usepackage{amsmath}
\usepackage{bbm}
\usepackage{amsfonts}
\usepackage{amsthm}
\newtheorem{rem1}{Remark}[section]
\newtheorem{lem1}{Lemma}[section]
\newtheorem{cor1}{Corollary}[section]

\newtheorem{thm1}{Theorem}[section]

\begin{document}

\title{Semi-linearity of the non-linear Fourier transform of the defocusing NLS equation}\author{T. Kappeler\footnote{Supported in part 
by the Swiss National Science Foundation}, B. Schaad\footnote{Supported in part by the Swiss National Science Foundation}, 
P. Topalov\footnote{Supported in part by NSF DMS-0901443}}

\maketitle

\begin{abstract} \noindent 
In this paper we prove that the non-linear Fourier transform of the defocusing NLS equation on the circle is linear up to terms which  
are one order smoother.

\medskip



\end{abstract}

\section{Introduction}
Consider the defocusing non-linear Schr\"odinger (dNLS) equation
\begin{align}\label{1} 
i \partial_tu=-\partial_x^2 u +2|u|^2u
\end{align}
on the circle $\mathbb{T}= \mathbb{R}/\mathbb{Z}.$ According to \cite{Bo0}, the initial value problem of \eqref{1} is globally in 
time well-posed on the Sobolev spaces $H^s\equiv H^s(\mathbb{T},\mathbb{C})$ with $s\geq 0$ where 
\[
H^s:= \big\{u= \sum_{n\in \mathbb{Z}}\hat u(n)e^{2\pi i n x}\,\big|\, \|u\|_{H^s} < \infty \big\}
\]
with
\[ 
\|u\|_{H^s}:= \Big(\sum_{n\in \mathbb{Z}}\langle n\rangle^{2s} |\hat u (n)|^2\Big)^{\frac{1}{2}} , 
\quad \hat u(n)= \int_0^1 u(x) e^{-2\pi inx} dx,
\]
and $\langle n\rangle:= \max(1,|n|)$. Equation \eqref{1} can be written as a Hamiltonian PDE with the real subspace 
\[
H^s_r:= \big\{(u,\bar{u})\,\big|\,u\in H^s\big\}\subseteq H^s\times H^s
\]
of the complex space $H^s_c:= H^s \times H^s$ as a phase space and Poisson bracket 
\[
\left\{F,G\right\}:=-i \int_0^1 \big(\partial_u F\partial_{\bar{u}}G - \partial_{\bar u} F\partial_u G\big)\,dx
\]
where $F, G: H^s_r\to \mathbb{C}$ are $C^1$-smooth functionals with sufficiently smooth 
$L^2$-gradients. The $L^2$-gradients $\partial_u F$ and $\partial_{\bar u}F$
are defined in a standard way in terms of the $L^2$-gradients of $F$ with respect to the real and imaginary
parts of $u$. The dNLS equation then takes the form 
\[
\partial_t u= - i \partial_{\bar u}\mathcal{H}_{\rm NLS}
\] 
where $\mathcal{H}_{\rm NLS}$ is the dNLS Hamiltonian 
\[
\mathcal{H}_{\rm NLS}:=\int_0^1 \big(\partial_x u \partial_x{\bar u}+ u^2 {\bar u}^2 \big)\,dx.
\]
According to \cite{GK}, the dNLS equation is an integrable PDE in the following strong sense: on $H^0_r$ there exists a 
canonical real analytic coordinate transformation $\Phi$  so that the initial value problem of \eqref{1},  when expressed in
these coordinates, can be explicitly solved by quadrature. The transformation $\Phi$ is referred to as {\em non-linear Fourier transform}
or {\em Birkhoff map} for the dNLS equation. It is defined on a complex neighborhood $W$ of $H^0_r$ in $H^0_c$  and takes values in
a complex neighborhood of the real subspace $\mathfrak{h}^0_r$ of $\mathfrak{h}^0_c$ where for any $s\in \mathbb{R}$,
\[
\mathfrak{h}^s_r:=\big\{(z,\bar z)\,\big|\, z\in \mathfrak{h}^s\big\}\subseteq\mathfrak{h}^s_c:=\mathfrak{h}^s\times\mathfrak{h}^s,
\]
and 
\[
\mathfrak{h}^s\equiv \mathfrak{h}^s(\mathbb{Z}, \mathbb{C}):= 
\big\{z=\big(z(n)\big)_{n\in \mathbb{Z}}\subseteq \mathbb{C}\,\big|\, \|z\|_s < \infty\big\} 
\]
with 
\[
\|z\|_s:=\Big(\sum_{n\in \mathbb{Z}}\langle n\rangle^{2s} |z(n)|^2\Big)^{\frac{1}{2}}.
\]
The Poisson structure on $\mathfrak{h}^s_r$ is defined by the condition that the coordinates 
$\big(\big(z(n)\big)_{n\in \mathbb{Z}},\big(\overline{z}(n)\big)_{n\in \mathbb{Z}}\big)$ in $\mathfrak{h}_r^0$ satisfy the 
canonical  relations  $\left\{ z(n), \overline{z}(n) \right\}= -i$ whereas all other brackets between coordinate functions vanish.
The Birkhoff map $\Phi$ has the property that when restricted to $H^N_r$, $N\geq 1$,  it takes values in $\mathfrak{h}^N_r$.

\noindent In order to state our results, it is convenient to denote by $\mathcal{F}$ the following Fourier transform,
\[
\mathcal{F}:H^0_c\to \mathfrak{h}^0_c, \quad
(\varphi_1, \varphi_2) \mapsto \left(\left(-\hat\varphi_1(-n)\right)_{n\in\mathbb{Z}},\left(-\hat\varphi_2(n)\right)_{n\in\mathbb{Z}}\right).
\]
Note that for $\varphi\in H^0_r,\, \bar\varphi_2= \varphi_1$ implying that $\overline{\hat\varphi}_2(n)= \hat\varphi_1(-n)$ for any 
$n\in\mathbb{Z}.$ Hence $\mathcal{F}$ maps $H^0_r$ into $\mathfrak{h}^0_r$.
Clearly, for any $s\geq 0,\, \mathcal{F}:H^s_r \to\mathfrak{h}^s_r$ is an isometry. 

\medskip

\noindent By definition, a (possibly non-linear) map between Banach spaces is said to be bounded if the image of any 
bounded subset is bounded. 

\begin{thm1}\label{Theorem1.1}
For any $N\in \mathbb{Z}_{\geq 1},$ the restriction of
$\Phi-\mathcal{F}$ to $ H_r^N$ takes values in $\mathfrak{h}_r^{N+1}$. The map $\Phi-\mathcal{F} : H_r^N \to \mathfrak{h}_r^{N+1}$ is real analytic and bounded.
\end{thm1} 
\noindent Theorem \ref{Theorem1.1} yields the following two corollaries:
\begin{cor1}\label{Corollary1.2} (i)
For any $N\in \mathbb{Z}_{\geq 1},$ the restriction of $\Phi^{-1}- \mathcal{F}^{-1}$ to $\mathfrak{h}^N_r$ takes values in $H^{N+1}_r$. The map  $\Phi^{-1}- \mathcal{F}^{-1}: \mathfrak{h}^N_r\to H^{N+1}_r$ is real analytic and bounded.
(ii) For any
 $s \in \mathbb{R}_{\geq 1},$ the restriction of 
 $\Phi$ to $H^s_r$ takes values in  $\mathfrak{h}^s_r$. 
The map $\Phi:H^s_r\to  \mathfrak{h}^s_r$ is a real analytic and bounded diffeomorphism as is its inverse.

\end{cor1}
\begin{cor1}\label{Corollary1.3}
For any $s\in \mathbb{R}_{\geq 1}, \, \Phi:H^s_r\to \mathfrak{h}^s_r$ and $\Phi^{-1}: \mathfrak{h}^s_r \to H^s_r$ are weakly 
continuous. 
\end{cor1}

\medskip

\noindent {\em Method of proof:} The proof of Theorem \ref{Theorem1.1} is based on a new formula for the components of the Birkhoff map
 presented in Theorem \ref{Theorem2.2} in Section 2. The main ingredients of the proof of Theorem \ref{Theorem1.1}, given in Section 3, 
are sharp asymptotic estimates of various spectral quantities of the Zakharov-Shabat operator established recently in \cite{KST2}. 
This operator appears in the Lax pair formulation of the dNLS equation. In Section 3 one also finds the proofs of Corollary \ref{Corollary1.2} 
and Corollary \ref{Corollary1.3}.

\medskip

\noindent {\em Related work:} In \cite{KuPi}, Kuksin and Piatnitski initiated a study of random perturbations with damping of the KdV 
equation on the circle. It then was further continued by Kuksin \cite{Ku}. The purpose of these investigations are to describe how the KdV 
action variables evolve under certain perturbations of the KdV equation. To this end, the perturbed equation is expressed in KdV Birkhoff 
coordinates, constructed in \cite{KdV&KAM}. Up to highest order it is a linear differential equation if the non-linear part of the Birkhoff map is 
{\em 1-smoothing} (or, equivalently,  a {\em semi-linear} map). In \cite{KuPe}, Kuksin and Perelman succeeded in showing that near the equilibrium point, 
the non-linear part of the Birkhoff map is indeed 1-smoothing and conjectured that this holds true globally in phase space. 
In \cite{KST1} the authors show that this is indeed the case, i.e., that the Birkhoff map is semi-linear, and  apply their result to obtain various 
new features of solutions of the KdV equation on the circle.
In this paper and the subsequent one \cite{KST3}, we establish such  results for the defocusing NLS equation, another important non-linear 
dispersive evolution equation.

\section{Birkhoff map}
In this section we review the Birkhoff map, constructed in \cite{GK}, state and prove the new formula mentioned in the introduction as well as 
the asymptotic estimates of quantities appearing in this formula. 
In \cite{GK}, Birkhoff coordinates $x(n),y(n),\,n\in \mathbb{Z}$, were constructed on $H^0_r$. These coordinates are real valued and
 satisfy $\{x(n),y(n)\}= -1$ for any $n\in \mathbb{Z}$ whereas all other brackets between coordinate functions vanish. 
For our purposes, in this paper it is convenient to use complex coordinates
\[
z_1(n):=\frac{x(n)-iy(n)}{\sqrt{2}}\quad \text{and} 
\quad z_2(n):= \frac{x(n)+ iy(n)}{\sqrt{2}} \,(= \overline{z}_1(n)).
\]
Then $\big\{z_1(n), \overline{z}_1(n)\big\}=-i$ whereas all other brackets between coordinate functions vanish and the action variables 
$I_n:= \big( x(n)^2+y(n)^2\big)/2$ can be expressed as 
\[
I_n= z_1(n)\overline{z}_1(n)\quad 
\forall n\in \mathbb{Z}.
\]
The result on the Birkhoff map in \cite{GK} (see the Overview as well as Theorem 20.2)
then reads as follows:
\begin{thm1}\label{Theorem2.1}
There exists an open neighborhood $W\subseteq H^0_c$ of $H^0_r$ and a real analytic map $\Phi: W\to \mathfrak{h}^0_c, \; (\varphi_1,\varphi_2)\mapsto (z_1, z_2)$ so that the following properties hold:
\begin{itemize}
\item[(B1)] $\Phi$ is canonical, i.e., preserves the Poisson brackets. 
\item[(B2)] The restriction of $\Phi$ to $H^N_r, \, N\in \mathbb{Z}_{\geq 0}$, gives rise to a map $\Phi: H^N_r\to \mathfrak{h}^N_r$ that is onto and bianalytic. 
\item[(B3)] $\Phi$ defines global Birkhoff coordinates for the dNLS equation on $H^1_r$. More precisely, 
on $\mathfrak{h}^1_r$, the dNLS Hamiltonian $\mathcal{H}_{\rm NLS}\circ \Phi^{-1}$ is a real analytic function of the actions $I_n=z_1(n)\overline{z}_1(n), \, n\in \mathbb{Z},$ alone.
\item[(B4)] The differential of $\Phi$ at $0$ is the Fourier transform, $d_0\Phi= \mathcal{F}$.
\end{itemize}  
\end{thm1}
\begin{rem1}\label{Remark2.3}
The claim (B4) of Theorem \ref{Theorem2.1} follows from \cite{GK} (Theorem 17.2). 
\end{rem1}
\noindent To state our new formula for the Birkhoff coordinates we first need to introduce some more notation and results from \cite{GK}.
For $\varphi\in H^0_c$, denote by $L(\varphi)$ the Zakharov-Shabat operator,   
\[
L(\varphi) = i \begin{pmatrix}
1&0\\0&-1
\end{pmatrix}\partial_x + \begin{pmatrix}
0& \varphi_1 \\
\varphi_2 &0
\end{pmatrix}
\]
and by $M(x,\lambda)\equiv M(x,\lambda,\varphi)$ the fundamental solution,  
\[
M(x,\lambda)= \begin{pmatrix}
m_1(x,\lambda)& m_2(x,\lambda) \\
m_3(x,\lambda) & m_4(x,\lambda)
\end{pmatrix},\quad M(0,\lambda)= \begin{pmatrix}
1&0\\0&-1
\end{pmatrix}. 
\]
It satisfies $LM= \lambda M$ for any $\lambda\in \mathbb{C}$ and for any $x\in \mathbb{R},\,M(x,\lambda) $  is an entire function of $\lambda$. 
We need to consider the Dirichlet spectrum of $L(\varphi)$ on $[0,1]$ as well as the periodic spectrum of $L(\varphi)$ on $[0,2].$ Both spectra are discrete and can be listed (with their algebraic multiplicities) as sequences
\[
\cdots \preceq \mu_{n-1}\preceq \mu_n\preceq \mu_{n+1} \preceq \cdots\quad\text{and}\quad\cdots \preceq\lambda^-_{n}\preceq \lambda^+_n\preceq \lambda^-_{n+1} \preceq\lambda_{n+1}^+ \preceq \cdots
\]
in lexicographic ordering $\preceq$ in such a way that $\mu_n,\lambda_n^\pm= n\pi +o(1)$ as $|n|\to \infty$. 
For complex numbers $a,b$ we write $a\preceq b$,  if $[\operatorname{Re}a< \operatorname{Re}b]$ or $[\operatorname{Re}a=\operatorname{Re}b\;\text{and}\,\operatorname{Im}a\leq\operatorname{Im}b]$. 
Furthermore denote by $\Delta(\lambda) \,[\delta(\lambda)]$ the trace [anti-trace] of $M(1,\lambda)$,
\[
\Delta(\lambda)=m_1(1,\lambda)+m_4(1,\lambda)\quad 
\delta(\lambda)=m_2(1,\lambda)+m_3(1,\lambda).
\] 
By Floquet theory, $(\lambda_n^\pm)_{n\in \mathbb{Z}}$ are the roots (with multiplicities) of $\Delta^2(\lambda)-4,$ which admits the product representation (\cite{GK}, Lemma 6.8), 
\[\Delta^2(\lambda)-4= -4 \prod_{k\in \mathbb{Z}}\frac{(\lambda_k^+-\lambda)(\lambda_k^-- \lambda)}{\pi_k^2}\]
where $\pi_k= k\pi$ for $k\neq 0$ and $\pi_0=1$. Let $\dot\Delta(\lambda)= \partial_\lambda \Delta(\lambda).$ Its zeros can be listed (with their multiplicities) as a sequence  $\cdots \preceq \dot\lambda_{n-1}\preceq \dot\lambda_n\preceq \dot\lambda_{n+1} \preceq \cdots$ in lexicographic ordering such that $\dot\lambda_n=n\pi +o(1)$ as $|n|\to \infty$.
The entire function $\dot\Delta(\lambda)$ also admits a product representation \cite{GK}, Lemma 6.5, 
\[
\dot \Delta(\lambda)= 2\prod_{k\in \mathbb{Z}}\frac{\dot \lambda_k-\lambda}{\pi_k}.
\]
It is shown in \cite{GK}, Section 12, that there exists an open neighborhood $W$ of $H^0_r$ in $H^0_c$ so that any element in $W$ admits a sequence of pairwise disjoint discs, $(U_n)_{n\in \mathbb{Z}}$ with the property that for any $n\in \mathbb{Z}$
\begin{align*}
\mu_n,\dot\lambda_n\in U_n, &\quad\text{and}\quad 
[\lambda^-_k,\lambda^+_k] = \{(1-t)\lambda_n^-+ t\lambda_n^+ |\, 0\leq t\leq 1\}\subseteq U_n
\\ c^{-1}|m-n|\leq &\operatorname{dist}(U_m,U_n) \leq c|m-n|\quad \forall m\neq n
\end{align*}
for some $c\geq 1.$ It can be shown that such sequences of discs, referred to as isolating neighborhoods, can be chosen locally uniformly on $W$ and so that for $|n|$ sufficiently large $U_n=D_n$ where 
$D_n= \{\lambda\in \mathbb{C}|\, |\lambda-n\pi| < \frac{\pi}{4}\}.$ 
 For $\varphi \in W,$ the action variables $I_n,\, n\in \mathbb{Z},$ are then defined as follows (\cite{GK}, Section 13) 
\[I_n = \frac{1}{\pi}\int_{\Gamma_n}\frac{\dot\Delta(\lambda)}{\sqrt[c]{\Delta^2(\lambda)-4}}d\lambda\]
where $\Gamma_n\subseteq U_n$ is a contour around $G_n$ and  $\sqrt[c]{\Delta^2(\lambda)-4}$ denotes the canonical root 
\begin{align*}
\sqrt[c]{\Delta^2(\lambda)-4}=2i \prod_{n\in \mathbb{Z}}\frac{\sqrt[s]{(\lambda_k^+-\lambda)(\lambda_k^-- \lambda)}}{\pi_k}
\end{align*}
with $\sqrt[s]{(\lambda_k^+-\lambda)(\lambda_k^-- \lambda)}$ being the standard root, defined on $\mathbb{C}\setminus  [\lambda^-_k,\lambda^+_k]$  by requiring that for $|\lambda|$ sufficiently large
\[
\sqrt[s]{(\lambda_k^+-\lambda)(\lambda_k^-- \lambda)}= 
(\tau_k-\lambda)\sqrt[+]{1- \Big(\frac{\gamma_k/2}{\tau_k-\lambda}\Big)^2}.
\]
Here $\sqrt[+]{x}$ denotes the root on $ \mathbb{C}\setminus(-\infty,0]$ with $\sqrt[+]{1}=1.$
In \cite{GK}, Theorem 3.3, it is shown that for any $\varphi\in W$ and $n\in \mathbb{Z}$ with $\gamma_n:= \lambda_n^+-\lambda_n^-\neq 0,\, 4I_n /\gamma_n^2$ is in $ \mathbb{C}\setminus(-\infty,0]$ and that the function $\xi_n=\sqrt[+] {4I_n / \gamma_n^2}$ extends analytically to $W$ and satisfies locally uniformly on $W$ the asymptotic estimate 
\[
\xi_n=1+\ell^2_n
\]
meaning $(\xi_n-1)_{n\in\mathbb{Z}}$ is a sequence in $\ell^2(\mathbb{Z},\mathbb{C})=\mathfrak{h}^0.$ Furthermore we recall that the angle variables are defined in terms of the entire functions $\psi_n(\lambda),$
\[
\psi_n(\lambda)= - \frac{2}{\pi_n} \prod_{k\neq n}\frac{\sigma^n_k-\lambda}{\pi_k}
\]
where the normalizing factor $-\frac{2}{\pi_n}$ and the zeros $\sigma_k^n,\,k\neq n,$ are chosen in such a way that 
\[ 
\frac{1}{2\pi}\int_{\Gamma_m} \frac{\psi_n(\lambda)}{\sqrt[c]{\Delta^2(\lambda)-4}}d\lambda=\delta_{nm}\quad \forall n,m\in \mathbb{Z}.
\]
The functions $\psi_n$ are used to define
\begin{align}\label{2.1}
\beta_{n,k}= \int_{\lambda_{k}^-}^{\mu_k }\frac{\psi_n(\lambda)}{\sqrt[*]{\Delta^2(\lambda)-4}}d\lambda,\quad \forall k\neq n
\end{align}
where the path of integration is required to be in $U_k\setminus G_k$ except possibly its endpoints, 
but otherwise arbitrary and the $*$-root $\sqrt[*]{\Delta^2(\lambda)-4}$ is defined by 
\[
\sqrt[*]{\Delta^2(\mu_k)-4}= \delta(\mu_k).
\]
Note that by \cite{GK}, Lemma 6.6, $\delta^2(\mu_k)=\Delta^2(\mu_k)-4.$
(As $\int_{\lambda_k^-}^{\lambda_k^\pm}\frac{\psi_n(\lambda)}{\sqrt{\Delta^2(\lambda)-4}}d\lambda=0,$ the choice of the sign of $\sqrt{\Delta^2(\lambda)-4}$ doesn't matter for $\mu_k\in \{\lambda_k^\pm\}$.)
For $k=n,\,\beta_{n,n}$ is defined also by \eqref{2.1}, but in this case only mod($2\pi i$) and for $\varphi\in W$ with $\gamma_n\neq 0.$ In \cite{GK}, Lemma 15.1, it is shown that $\beta_{n,k}$, $k \ne n$, is analytic and satisfies
\begin{align}
\label{2.1bis} \beta_{n,k}= O\left(\frac{|\tau_k-\mu_k| + |\gamma_k|}{|n-k|}\right)
\end{align} 
locally uniformly on $W,$ implying that $\beta_n= \sum_{k\neq n}\beta_{n,k}$ is absolutely summable due to the asymptotic estimates of $\tau_k-\mu_k$ and $\gamma_k.$ The functions $\beta_n$ are shown to be analytic on $W$  and to satisfy $\beta_n=o(1),\, |n|\to \infty,$ locally uniformly on $W$ in \cite{GK}, Theorem 15.3. 
In \cite{GK}, Section 16, the Birkhoff coordinates are then defined as
\begin{align}
\label{2.2} z_1(n)= \xi_ne^{-i\beta_n} z_n^-, \quad z_2(n)= \xi_n e^{i\beta_n} z_n^+, \quad n\in \mathbb{Z}
\end{align}
where $z_n^{\pm}$ are defined for $\varphi\in W$ with $\gamma_n \neq 0$ by $z_n^{\pm}= \frac{\gamma_n}{2}e^{\pm i \beta_{n,n}}$ and then shown to extend as analytic functions to $W$ in \cite{GK}, Proposition 16.5. (The factor $1/2$ in the definition of $z_n^{\pm}$  is convenient for our setting.)
Write $\delta(\mu_n)=\sqrt[*]{\Delta^2(\mu_n)-4}$ as a product
\begin{align}
\label{2.3} \delta(\mu_n)= 2i \sqrt[*]{(\lambda_n^+-\mu_n)(\lambda_n^--\mu_n)} \, \delta_n(\mu_n)
\end{align}
with the sign of $\sqrt[*]{(\lambda_n^+-\mu_n)(\lambda_n^--\mu_n)}$ defined by the latter equation and $\delta_n(\mu_n)$ given by \begin{align}
\label{2.3bis} \delta_n(\mu_n)= \frac{1}{\pi_n} \prod_{k\neq n} \frac{\sqrt[s]{(\lambda_k^+ -\mu_n)(\lambda_k^- -\mu_n)}}{\pi_k}.
\end{align}
By \cite{GK}, Lemma 12.7, $\delta_n(\mu_n)$ is analytic on $W$. Similarly write 
\[\frac{\psi_n(\lambda)}{\sqrt[*]{\Delta^2(\lambda)-4}}= \frac{i}{\sqrt[*]{(\lambda_n^+-\lambda)(\lambda_n^--\lambda)}}\zeta_n(\lambda)\]
with \begin{align}
\label{2.4} \zeta_n(\lambda):= \prod_{k\neq n} \frac{\sigma_k^n-\lambda}{\sqrt[s]{(\lambda_k^+-\lambda)(\lambda_k^--\lambda)}}
.\end{align}
By \cite{GK}, Lemma 12.10, and the asymptotics of $\mu_n,$ \begin{align}
\label{2.5} \delta_n(\mu_n)^{-1}= (-1)^n + \ell^2_n
\end{align} 
and by \cite{GK}, Lemma 12.2, as $|n|\to \infty$
\begin{align}\label{2.6}
\zeta_n(\lambda)= 1 + \ell_n^2 \quad \text{uniformly for}\; \lambda\in U_n \,.
\end{align}
Both asymptotic estimates, hold locally uniformly on $W$. 
With these notations we can now state the new formulas for $z_n^\pm:$

\begin{thm1}\label{Theorem2.2}
For any $\varphi\in W$ 
and $n\in \mathbb{Z}$\[z_n^\pm = \left((\tau_n-\mu_n)\mp i\frac{\delta(\mu_n)}{2\delta_n(\mu_n)}
\right)\eta_n^\pm\]
where $\eta_n^\pm$ are defined by 
\[\eta_n^\pm=\begin{cases}
\exp\left(\mp \int_{\lambda_n^-}^{\mu_n} \frac{\zeta_n(\lambda)-1}{\sqrt[*]{(\lambda_n^+-\lambda)(\lambda_n^--\lambda)}}d\lambda \right)\quad &\text{if} \; \gamma_n\neq 0 \\
\exp\left(\mp \varepsilon_n\int_{\tau_n}^{\mu_n} \frac{\zeta_n(\lambda)-\zeta_n(\tau_n)}{\tau_n-\lambda}d\lambda \right)\quad&\text{if} \; \gamma_n= 0,\; \mu_n\neq \tau_n \\
1  &\text{if} \; \gamma_n=0,\; \mu_n=\tau_n
\end{cases}
\]
with $\varepsilon_n\in \{\pm 1\}$ given by $\varepsilon_n= \frac{\delta(\mu_n)}{\sqrt[c]{\Delta^2(\mu_n)-4}}$ if $\gamma_n=0,\mu_n\neq \tau_n.$ The functions $\eta_n^\pm$ are analytic and satisfy
\begin{align}
\label{2.7} \eta_n^\pm= 1+ O\big( |\gamma_n|+ |\mu_n-\tau_n|\big)\quad\text{as}\; |n| \to \infty
\end{align}
locally uniformly on $W.$
\end{thm1}
Expressed in an informal way, Theorem \ref{Theorem2.2}, when combined with \eqref{2.1bis} and \eqref{2.2},  says that
\begin{align}
\label{2.8} z_1(n)\sim (\tau_n-\mu_n)+ i(-1)^n\frac{\delta(\mu_n)}{2} \\\label{2.9} z_2(n)\sim (\tau_n-\mu_n)- i(-1)^n\frac{\delta(\mu_n)}{2}
\end{align}
up to error terms of order $o(1)$ as $|n|\to \infty.$ Note that $\tau_n,\mu_n,\delta(\mu_n)$ and $\delta_n(\mu_n)$ are analytic on $W.$
\begin{proof}
[Proof of Theorem 2.2]
For any $n\in \mathbb{Z},$ let $Z_n:= \{\varphi \in W|\, \gamma_n=0\}. $ Recall that for $\varphi\in W\setminus Z_n,$ $z_n^-= \frac{\gamma_n}{2}e^{-i\beta_{n,n}}$ with
\[
\beta_{n,n}= \int_{\lambda_{n}^-}^{\mu_n }\frac{\psi_n(\lambda)}{\sqrt[*]{\Delta^2(\lambda)-4}}d\lambda \quad\text{mod} (2\pi i) \, .
\]
(As $-i\int_{\lambda_{n}^-}^{\lambda_n^+ }\frac{\psi_n(\lambda)}{\sqrt{\Delta(\lambda)^2-4}}d\lambda=\pi\; (\text{mod} 2\pi)$   independently of the choice of the sign of the root,
the sign of the root in the integrand is irrelevant if $\mu_n\in \{\lambda_n^\pm\}.$ )
With $\zeta_n(\lambda)$ given by \eqref{2.4}, one has 
\[-i\beta_{n,n}= 
\int_{\lambda_{n}^-}^{\mu_n }\frac{d\lambda}{\sqrt[*]{(\lambda_n^+-\lambda)(\lambda_n^--\lambda)}}+ 
\int_{\lambda_n^-}^{\mu_n} \frac{\zeta_n(\lambda)-1}{\sqrt[*]{(\lambda_n^+-\lambda)(\lambda_n^--\lambda)}}d\lambda.\]
We claim that for $\varepsilon\in \{1,-1\}$ and $\mu\in \mathbb{C}\setminus G_n$
\[
\frac{\gamma_n}{2}\exp\left(\varepsilon \int_{\lambda_{n}^-}^{\mu }
\frac{d\lambda}{\sqrt[*]{(\lambda_n^+-\lambda)(\lambda_n^--\lambda)}}\right)= 
(\tau_n-\mu) - \varepsilon\sqrt[*]{(\lambda_n^+-\mu)(\lambda_n^--\mu)}. 
\]
Indeed,  both sides are analytic univalent functions on $ \mathbb{C}\setminus G_n$ with limit $\frac{\gamma_n}{2}$ as $\mu \to \lambda_n^-,$ satisfying the differential equation 
\[f'(\mu)= \frac{\varepsilon}{\sqrt[*]{(\lambda_n^+-\mu)(\lambda_n^--\mu)}}f(\mu)\quad \forall \mu\in \mathbb{C}\setminus G_n.\]
When combined with \eqref{2.3} one then gets 
\[z_n^-=\frac{\gamma_n}{2}e^{-i \beta_{n,n}}= \left((\tau_n-\mu_n)+ i\frac{\delta(\mu_n)}{2\delta_n(\mu_n)}\right)\eta_n^-(\mu_n).\]
To see that $\eta_n^-$ is analytic on $W$ note first that for any $\varphi \in W\setminus Z_n,$ by the definition $\psi_n,$
\[\int_{\lambda_{n}^-}^{\lambda_n^+ }\frac{\zeta_n(\lambda)}{\sqrt{(\lambda_n^+-\lambda)(\lambda_n^--\lambda)}}d\lambda=\pi \;(\text{mod} 2\pi)\]
for any choice of the sign of the roots. As
\[\int_{\lambda_{n}^-}^{\lambda_n^+ }\frac{1}{\sqrt{(\lambda_n^+-\lambda)(\lambda_n^--\lambda)}}d\lambda=\pi \;(\text{mod} 2\pi)\]
by a straightforward computation it follows that 
\[
\exp\left( \int_{\lambda_n^-}^{\lambda_n^+} 
\frac{\zeta_n(\lambda)-1}{\sqrt[*]{(\lambda_n^+-\lambda)(\lambda_n^--\lambda)}}\,d\lambda\right)=1
\]
and hence for any $\varphi \in W\setminus Z_n$ one has 
\begin{align}\label{2.10}
\exp\left(\int_{\lambda_n^-}^{\mu_n} 
\frac{\zeta_n(\lambda)-1}{\sqrt[*]{(\lambda_n^+-\lambda)(\lambda_n^--\lambda)}}\,d\lambda\right)
=\exp\left(\int_{\lambda_n^+}^{\mu_n} 
\frac{\zeta_n(\lambda)-1}{\sqrt[*]{(\lambda_n^+-\lambda)(\lambda_n^--\lambda)}}\,d\lambda\right).
\end{align}
Arguing as in the proof of \cite{GK}, Lemma 15.2, one sees that $\eta_n^-$ is analytic on $W\setminus Z_n$. To show that it is analytic on all of $W$ one argues as in the proof of \cite{GK}, Proposition 16.5: one verifies that $\eta_n^-$ is continuous on $W$ and that its restriction to $Z_n$ is weakly analytic to conclude from \cite{GK}, Theorem A.6, that $\eta_n^-$ is analytic on $W.$ The details are left to the reader.

\noindent To prove the claimed asymptotics for $\eta_n^-$ we first consider the case where $\varphi\in Z_n.$ In the case $\mu_n= \tau_n$ one has $\eta_n^-=1.$ If $\mu_n\neq \tau_n,$ then with the parametrization $\lambda(t)= \tau_n +t(\mu_n-\tau_n),\; 0\leq t\leq 1,$ one gets
\[
\left|\int_{\tau_n}^{\mu_n}\frac{\zeta_n(\lambda)-\zeta_n(\tau_n)}{\tau_n-\lambda}d\lambda \right| \leq \int_0^1\left|\frac{\zeta_n(\lambda)-\zeta_n(\tau_n)}{\tau_n-\lambda} \right||\mu_n-\tau_n|dt
\]
yielding the estimate 
\[
\eta_n^-= \exp \left(O(|\mu_n-\tau_n|)\right)= 1+O(|\mu_n-\tau_n|).
\]
In the case $\varphi \in W\setminus Z_n,$ write $\eta_n^-= \exp(I)\exp(II)$
where 
\[
I=\int_{\lambda_n^-}^{\mu_n}\frac{\zeta_n(\lambda)-\zeta_n(\tau_n)}{\tau_n-\lambda}\frac{\tau_n-\lambda}{\sqrt[*]{(\lambda_n^+-\lambda)(\lambda_n^--\lambda)}}d\lambda 
\]
and 
\[
II= \int_{\lambda_n^-}^{\mu_n}\frac{\zeta_n(\tau_n)-1}{\sqrt[*]{(\lambda_n^+-\lambda)(\lambda_n^--\lambda)}}d\lambda.
\]
Towards $I,$ with the path of integration given by 
$\lambda(t)= \lambda_n^- +t(\mu_n-\lambda_n^-)$ one has
\[
|I|\leq \int_0^1 \left|\frac{\zeta_n(\lambda)-\zeta_n(\tau_n)}{\tau_n-\lambda}\right|\left|\frac{\tau_n-\lambda}{\lambda_n^+-\lambda}\right|^{\frac{1}{2}}\left|\frac{\tau_n-\lambda}{\lambda_n^--\lambda}\right|^{\frac{1}{2}}\left|\mu_n-\lambda_n^-\right|dt.
\]
By \eqref{2.10} we can assume without loss of generality that $|\lambda_n^--\mu_n|\leq |\lambda_n^+-\mu_n|, $ as otherwise we can switch the roles of $\lambda_n^-$ and $\lambda_n^+.$ Hence we can assume that 
\[
|\lambda_n^+-\lambda(t)|\geq \left| \gamma_n / 2\right|\quad \forall 0\leq t\leq 1
\]
and thus
\[
\Big|\frac{\tau_n-\lambda(t)}{\lambda_n^+-\lambda(t)}\Big|^{\frac{1}{2}} \leq 
\Big(1+ \Big|\frac{\gamma_n / 2}{\lambda_n^+-\lambda(t)}\Big|^{\frac{1}{2}}\Big) \leq \sqrt{2}.
\] 
Furthermore, $\big|\frac{\tau_n-\lambda(t)}{\lambda_n^--\lambda(t)}\big|^{\frac{1}{2}} \leq \frac{\left(|\frac{\gamma_n}{2}|+ |\mu_n-\lambda_n^-|\right)^{\frac{1}{2}}}{\sqrt{t}|\mu_n-\lambda_n^-|^{\frac{1}{2}}}$. 
Hence 
\[
|I|\leq \Big(\int_0^1\left|\frac{\zeta_n(\lambda)-\zeta_n(\tau_n)}{\tau_n-\lambda}\right|\frac{dt}{\sqrt[+]{t}} \Big)\sqrt{2}
\Big(\Big|\frac{\gamma_n}{2}\Big|+ |\mu_n-\lambda_n^-|^{\frac{1}{2}}\Big) |\mu_n-\lambda_n^-|^{\frac{1}{2}}.
\]
As $\zeta_n(\lambda)$ is analytic in $\lambda \in U_n$ (cf \cite{GK}, Corollary 12.8) one has $\frac{\zeta_n(\lambda)-\zeta_n(\tau_n)}{\tau_n-\lambda}= O(1)$ yielding altogether the estimate 
\[
I= O\left(|\gamma_n|^{\frac{1}{2}}|\mu_n-\lambda_n^-|^{\frac{1}{2}}+ |\mu_n-\lambda_n^-|\right).
\]

\noindent Towards $II$, note that by \cite{GK}, Corollary 16.3, 
$\, \zeta_n(\tau_n)-1= O(\gamma_n).$ Again assuming that 
$|\lambda_n^+-\mu_n|\geq \left|\gamma_n / 2\right|$ and using the path of integration $\lambda(t)= \lambda_n^- +t(\mu_n-\lambda_n^-),\; 0\leq t\leq 1,$ one gets 
\[
|II|= C|\gamma_n|\int_0^1 \frac{1}{\left|\gamma_n / 2 \right|}\frac{1}{\sqrt[+]{t}|\mu_n-\lambda_n^-|^{\frac{1}{2}}}|\mu_n-\lambda_n^-|^{\frac{1}{2}} dt
\]
implying that $II= O\left(|\gamma_n|^{\frac{1}{2}}|\mu_n-\lambda_n^-|^{\frac{1}{2}}\right).$ Combining the estimates of $I$ and $II$ then yields 
\[
\eta_n^-= \exp\big(|\gamma_n|+|\mu_n-\tau_n|\big)= 1+O\big(|\gamma_n|+|\mu_n-\tau_n|\big).
\]
Going through the arguments of the proof of these estimates one sees that $\eta_n^-= 1+ O\big(|\gamma_n|+|\mu_n-\tau_n|\big)$ locally 
uniformly on $W.$ The proof of the claimed results for $\eta_n^+$ is of course similar.  
\end{proof}

\section{Proofs of the main results}
In this section we prove Theorem \ref{Theorem1.1}, Corollary \ref{Corollary1.2}, and Corollary \ref{Corollary1.3}. Without further  reference we use the notation introduced in the previous sections. The main ingredient in the proof of Theorem \ref{Theorem1.1} are the following asymptotic estimates, proved in \cite{KST2}: 
\begin{thm1}\label{Theorem3.1}
For $\varphi$ in $W\cap H^N_c$ with $N\in \mathbb{Z}_{\geq 1},$ 
\begin{itemize}
\item[(i)]  $\tau_n-\mu_n= - \left(\hat\varphi_1(-n)+\hat\varphi_2(n)\right)/2 + \frac{1}{n^{N+1}}\ell^2_n$
\item[(ii)] $\delta(\mu_n)= (-1)^n i \left(\hat\varphi_1(-n)-\hat\varphi_2(n)\right)+\frac{1}{n^{N+1}}\ell^2_n $
\item[(iii)] $\gamma_n= \frac{1}{n^{N}}\ell^2_n $
\end{itemize}
locally uniformly and uniformly on bounded subsets of $H^N_r.$
\end{thm1}
In the following two lemmas we establish a few additional asymptotic estimates for $\varphi \in W\cap H^1_c,$ needed for the proof of Theorem \ref{Theorem1.1}.
\begin{lem1}\label{Lemma3.2} On $W\cap H^1_c$ (i) $\beta_n=O\left( \frac{1}{n}  \right), $ 
(ii) $\eta_n^\pm=1+ O\left( \frac{1}{n}  \right),$ and 
(iii) $\xi_n^\pm=1+ O\left( \frac{1}{n}  \right)$ locally uniformly and uniformly on bounded subsets of $H^1_r.$
\end{lem1}
\begin{proof}
Recall that $\beta_n= \sum_{k\neq n}\beta_{n,k}$ with $\beta_{n,k}$ given by \eqref{2.1}
and satisfying $\beta_{n,k}= O\big(\frac{|\gamma_k|+|\mu_k-\tau_k|}{|k-n|}\big)$ 
by \eqref{2.1bis}. By Theorem \ref{Theorem3.1} it then follows that $\beta_{n,k}= O\big(\frac{1}{k}\frac{1}{k-n}\ell^2_k
\big),$ implying that $\beta_n=O\left( \frac{1}{n}  \right). $ 
Similarly, by Theorem \ref{Theorem2.2} and Theorem \ref{Theorem3.1} it follows that $\eta_n^\pm=1+ O\left( \frac{1}{n}  \right).$ Going through the arguments of the proof  and taking into account that the Sobolev embedding $H^1_c \hookrightarrow H_c^0$ is compact it follows that these estimates hold locally uniformly on $W\cap H^1_c$ and uniformly on bounded subsets of $H^1_r.$ Hence items (i) and (ii) are proved. 
Towards (iii), introduce for $\varphi \in W$ and $n\in \mathbb{Z},\; \chi_n(\lambda)= \prod_{k\neq n} \frac{\dot\lambda_k-\lambda}{\sqrt[s]{(\lambda_k^+-\lambda)(\lambda_k^--\lambda)}}$ and if $\gamma_n\neq 0, \; t_n := \frac{\dot\lambda_n-\tau_n}{\gamma_n/2}.$
By \cite{GK}, Theorem 13.3, 
\begin{align*}
\xi_n^2= \begin{cases}
\frac{2}{\pi} \int_{-1}^1\frac{(t-t_n)^2}{\sqrt[+]{1-t^2}}\chi_n(\tau_n +t\frac{\gamma_n}{2})d t & \text{if}\; \gamma_n\neq 0 \\
\chi_n(\tau_n) & \text{if}\; \gamma_n= 0 
\end{cases}
\end{align*}
If $\gamma_n\neq 0,$ write $\xi_n^2$ as a sum $I+II$ where 
\[I=\frac{2}{\pi} \int_{-1}^1\frac{(t-t_n)^2}{\sqrt[+]{1-t^2}}d t,\quad II= \frac{2}{\pi} \int_{-1}^1\frac{(t-t_n)^2}{\sqrt[+]{1-t^2}}\left(\chi_n(\tau_n +t\frac{\gamma_n}{2})-1\right)d t.\]
Let us first analyze $I.$ One computes 
\[
I= t_n^2 \frac{2}{\pi} \int_{-1}^1\frac{dt}{\sqrt[+]{1-t^2}}+ \frac{2}{\pi} \int_{-1}^1\frac{t^2}{\sqrt[+]{1-t^2}}d t= 2 t_n^2 +1.
\]
By \cite{GK}, Lemma 6.9, $t_n= O(\gamma_n)$ and thus 
\[I =1+ O(\gamma_n^2).\]
Towards term $II$, we claim that on $W\cap H^1_c$
\begin{align}\label{3.1} 
\chi_n\Big(\tau_n +t\frac{\gamma_n}{2}\Big)=1+ O\left(\frac{1}{n}\right)
\end{align} 
uniformly for $0\leq t\leq 1.$ 
To show this, choose $K\geq 1$ so that for any $|k|\geq K$ with $k\neq n,$ 
\[
\left| \frac{\gamma_k / 2}{\tau_k-\lambda}
\right| \leq \frac{1}{2}\quad \forall \lambda\in U_n.
\] 
Then for any $\lambda \in U_n,\; \sqrt[s]{(\lambda_k^+-\lambda)(\lambda_k^--\lambda)}= (\tau_k-\lambda)
\sqrt[+]{1- \big(\frac{\gamma_k / 2}{\tau_k-\lambda}\big)^2}$ and 
\[ 
\frac{\dot\lambda_k-\lambda}{\sqrt[s]{(\lambda_k^+-\lambda)(\lambda_k^--\lambda)}}=
\frac{\dot\lambda_k-\lambda}{\tau_k-\lambda}\Big(1- \Big(\frac{\gamma_k / 2}{\tau_k-\lambda}\Big)^2\Big)^{-\frac{1}{2}}
\]
yielding 
\begin{align*}
&\prod_{|k|\geq K, k\neq n}\frac{\dot\lambda_k-\lambda}{\sqrt[s]{(\lambda_k^+-\lambda)(\lambda_k^--\lambda)}}\\= 
\, &\left(\prod_{|k|\geq K, k\neq n}\Big(1+\frac{\dot\lambda_k-\tau_k}{\tau_k-\lambda}\Big)\right)
\cdot \left(\prod_{|k|\geq K, k\neq n}
\Big(1- \Big(\frac{\gamma_k / 2}{\tau_k-\lambda}\Big)^2\Big)\right)^{-\frac{1}{2}} 
\end{align*}
By \cite{GK}, Lemma C.2, uniformly for $\lambda\in U_n,$
\begin{align}\nonumber
\Big|\prod_{|k|\geq K, k\neq n}\Big(1+\frac{\dot\lambda_k-\tau_k}{\tau_k-\lambda}\Big)-1\Big|
\leq &\,C \sum_{|k|\geq K, k\neq n}\frac{|\dot\lambda_k-\tau_k|}{|k-n|}\\\label{3.2} 
\leq&\,C \sum_{|k|\geq K, k\neq n}\frac{|\gamma_k|^2}{|k-n|}
\end{align} 
where for the latter inequality we used that 
$\dot\lambda_k-\tau_k=O(\gamma_k^2)$ by \cite{GK}, Lemma 6.9. Again using \cite{GK}, Lemma C.2, one also has 
\begin{align}\label{3.3}
\Big|\prod_{|k|\geq K, k\neq n}\Big(1-\Big(\frac{\gamma_k / 2}{\tau_k-\lambda}\Big)^2\Big)-1 \Big|
\leq C \sum_{|k|\geq K, k\neq n}\Big|\frac{\gamma_k}{k-n}\Big|^2
\end{align} 
uniformly for $\lambda \in U_n.$ Finally, as for any $k\neq n$ 
\begin{align}\label{3.4}
\frac{\dot\lambda_k-\lambda}{\sqrt[s]{(\lambda_k^+-\lambda)(\lambda_k^--\lambda)}}= 1+ O\Big(\frac{1}{n-k}\Big)
= 1+ O\Big(\frac{1}{n}\Big)
\end{align}
uniformly for $\lambda \in U_n,$ the finite product $\prod_{|k|< K, k\neq n}\frac{\dot\lambda_k-\lambda}{\sqrt[s]{(\lambda_k^+-\lambda)(\lambda_k^--\lambda)}}$  is $1+O\big(\frac{1}{n}\big) $ uniformly for $\lambda \in U_n.$ Altogether it then follows from the asymptotics of $\gamma_n$ of Theorem \ref{Theorem3.1}, \eqref{3.1} holds uniformly for $0\leq t\leq 1,$ locally uniformly on $W\cap H^1_c$ and uniformly on bounded subsets of $H^1_r.$ As a consequence 
$II= 1+ O(\frac{1}{n})$ as well as $\chi_n(\tau_n)=1+O(\frac{1}{n}).$ 
In all we have shown that $\xi_n^2= 1+O\left(\frac{1}{n}\right)$ and thus $\xi_n= 1+O\left(\frac{1}{n}\right).$ 
Going through the arguments of the proof one verifies that the estimate holds locally uniformly on $W\cap H^1_c$ 
and uniformly on bounded subsets of $H^1_r.$
\end{proof}
\begin{lem1}
\label{Lemma3.3} On $W\cap H^1_c,\; \delta_n(\mu_n)^{-1}= (-1)^n+ O\left(\frac{1}{n}\right)$ locally uniformly and uniformly on 
bounded subsets of $H^1_r.$  
\end{lem1}
\begin{proof}
By \eqref{2.5}, $\delta_n(\mu_n)^{-1}= (-1)^n+ \ell_n^2.$ As by Theorem \ref{Theorem3.1}, $\mu_n= n\pi + O\left(\frac{1}{n}\right),$ 
\[
\frac{\sin \mu_n}{\mu_n-n\pi}= (-1)^n\frac{\sin (\mu_n-n\pi)}{\mu_n-n\pi} =(-1)^n+O\left((\mu_n-n\pi)^2\right) 
=  (-1)^n+O\Big(\frac{1}{n^2}\Big),
\]
yielding 
$\frac{\sin \mu_n}{\mu_n-n\pi}\frac{1}{\delta_n(\mu_n)}=  (-1)^n\frac{1}{\delta_n(\mu_n)}+O\big(\frac{1}{n^2}\big)$.
It therefore is to show that 
\begin{align*}
\frac{\sin \mu_n}{\mu_n-n\pi}\frac{1}{\delta_n(\mu_n)}
= 1+O\Big(\frac{1}{n}\Big).
\end{align*}
Using the product representations 
\[
\frac{1}{\delta_n(\mu_n)}= \pi_n \prod_{k\neq n} \frac{\pi_k}{\sqrt[s]{(\lambda_k^+-\mu_n)(\lambda_k^--\mu_n)}},
\qquad 
\frac{\sin \mu_n}{\mu_n-n\pi}= \frac{1}{\pi_n}\prod_{k\neq n} \frac{k\pi-\mu_n}{\pi_k}
\]
one sees that
\begin{align*}
\frac{\sin \mu_n}{\mu_n-n\pi}\frac{1}{\delta_n(\mu_n)}
= \prod_{k\neq n} \frac{k\pi-\mu_n}{\sqrt[s]{(\lambda_k^+-\mu_n)(\lambda_k^--\mu_n)}}.
\end{align*}
As in the proof of item (iii) of Lemma \ref{Lemma3.2}, choose $K\geq 1$ so that for any $|k|\geq K$ with 
$k\neq n,\; \frac{|\gamma_k /2 |}{\tau_k-\mu_n|} \leq \frac{1}{2}$ so that
\[
\sqrt[s]{(\lambda_k^+-\mu_n)(\lambda_k^-- \mu_n)}= 
(\tau_k-\mu_n)\sqrt[+]{1- \left(\frac{\gamma_k / 2}{\tau_k-\mu_n}\right)^2}
\]
implying that $\prod_{|k|\geq K, k\neq n} \frac{k\pi-\mu_n}{\sqrt[s]{(\lambda_k^+-\mu_n)(\lambda_k^--\mu_n)}}$ equals
\[
\left(\prod_{|k|\geq K, k\neq n}\Big(1-\frac{\tau_k-k\pi}{\tau_k-\mu_n}\Big)\right)
\left(\prod_{|k|\geq K, k\neq n}
\Big(1- \Big(\frac{\gamma_k / 2}{\tau_k-\mu_n}\Big)^2\Big)\right)^{-\frac{1}{2}}\,.
 \]
By \eqref{3.3}, 
\begin{align}\label{3.5}
\Big|\prod_{|k|\geq K, k\neq n}\Big(1- \Big(\frac{\gamma_k/2}{\tau_k-\mu_n}\Big)^2\Big)-1\Big| 
\leq C \sum_{|k|\geq K, k\neq n}\Big|\frac{\gamma_k}{k-n}\Big|^2
\end{align}
 and by \eqref{3.4} 
 \begin{align}
 \label{3.6}\prod_{|k|< K, k\neq n} \frac{k\pi-\mu_n}{\sqrt[s]{(\lambda_k^+-\mu_n)(\lambda_k^--\mu_n)}}= 1+ O\Big(\frac{1}{n}\Big).
 \end{align}
 The estimate of $\prod_{|k|\geq K, k\neq n}\big(1-\frac{\tau_k-k\pi}{\tau_k-\mu_n}\big)$ is more delicate. Due to the asymptotics of $\tau_k$ on $H^1_c$ (cf \cite{KST2})
 \begin{align}
 \label{3.8} \tau_k= k\pi +\frac{c_1}{k}+ O\Big(\frac{1}{k^2}\Big),\qquad c_1= \frac{1}{2\pi}\int_0^1\varphi_1(x)\varphi_2(x)dx\, .
 \end{align}
 Hence $|\frac{\tau_k-k\pi}{\tau_k-\mu_n}|\leq C\frac{1}{|k|}\frac{1}{|k-n|},\; \forall k\neq n,0.$ Thus $\big(\frac{\tau_k-k\pi}{\tau_k-\mu_n}\big)_{k\in \mathbb{Z}}$ is an $\ell^1$--sequence and, by \cite{GK}, Lemma C.2, it follows that
 $ P:=\big|\prod_{|k|\geq K, k\neq n}\big(1-\frac{\tau_k-k\pi}{\tau_k-\mu_n}\big)-1\big| $ satisfies
 \begin{align}\label{3.9} 
P \, \leq   C \,  \Big|\sum_{|k|\geq K, k\neq n}\frac{\tau_k-k\pi}{\tau_k-\mu_n}\Big| +
C \sum_{|k| \geq K, k\neq n}\Big|\frac{\tau_k-k\pi}{\tau_k-\mu_n}\Big|^2.
\end{align}  
Clearly 
\[ 
\sum_{|k|\geq K, k\neq n}\Big|\frac{\tau_k-k\pi}{\tau_k-\mu_n}\Big|^2
\leq C\sum_{|k|\geq K, k\neq n} \frac{1}{k^2}\frac{1}{|k-n|^2}
= O\Big(\frac{1}{n^2}\Big)
\]
whereas due to \eqref{3.8},
\[
\sum_{|k|\geq K, k\neq n}\frac{\tau_k-k\pi}{\tau_k-\mu_n}
= c_1 \sum_{|k|\geq K, k\neq n}\frac{1}{k}\frac{1}{\tau_k-\mu_n} 
+O\Big(\frac{1}{n}\Big).
\]
Furthermore, 
\begin{align*}
\frac{1}{\tau_k-\mu_n} = &\frac{1}{k\pi-n\pi} + \frac{(\mu_n-n\pi)+ (\tau_k-k\pi)}{(\tau_k-\mu_n)(k\pi-n\pi)}\\
= & \frac{1}{k\pi-n\pi} +O\Big(\Big(\frac{1}{|n|}+\frac{1}{|k|}\Big)\frac{1}{|k-n|^2}\Big).
\end{align*}
The asymtotics $\mu_n-n\pi=O\left(\frac{1}{n}\right)$ on $W\cap H^1_c$ (\cite{KST2}) then imply that 
\[ 
c_1 \sum_{|k|\geq K, k\neq n}\frac{1}{k}\frac{1}{\tau_k-\mu_n}
= \frac{c_1}{\pi} \sum_{|k|\geq K, k\neq n}\frac{1}{k(k-n)}+ O\Big(\frac{1}{n^2}\Big).
\]
As $\sum_{0<|k|< K, k\neq n}\frac{1}{k(k-n)}=O\left(\frac{1}{n}\right)$ it remains to prove that $\sum_{ k\neq n,0}\frac{1}{k(k-n)}=O\left(\frac{1}{n}\right)$ which is a classical result.
Substituting the estimates obtained into \eqref{3.9} lead to 
\begin{align}
\label{3.10}
\prod_{|k|\geq K, k\neq n}\Big(1-\frac{\tau_k-k\pi}{\tau_k-\mu_n}\Big)
= 1+ O\Big(\frac{1}{n}\Big).
\end{align}
Combining \eqref{3.5}, \eqref{3.6} and \eqref{3.10} then shows that 
\[
\prod_{ k\neq n} \frac{k\pi-\mu_n}{\sqrt[s]{(\lambda_k^+-\mu_n)(\lambda_k^--\mu_n)}}
= 1+ O\Big(\frac{1}{n}\Big).
\]
Going through the arguments of the proof one verifies that this estimate holds locally uniformly on $W\cap H^1_c$ and uniformly on 
bounded subsets of $H^1_r$.
\end{proof}
\begin{proof}
[Proof of Theorem \ref{Theorem1.1}] By \eqref{2.2} and Theorem \ref{Theorem2.2} one has for any $\varphi\in W,$
\[
z_1(n)= \xi_ne^{-i\beta_n}\Big((\tau_n-\mu_n)+ i\frac{\delta(\mu_n)}{2\delta_n(\mu_n)}\Big)\eta_n^-.
\]
By Lemma \ref{Lemma3.2}, for any $\varphi \in W\cap H^1_c$, one has that $\xi_n,\, e^{-i\beta_n},$ and $\eta_n^-$ satisfy the estimate $1+ O\big(\frac{1}{n}\big)$ whereas by Lemma \ref{Lemma3.3}, $\delta_n(\mu_n)^{-1}= (-1)^n+ O\big(\frac{1}{n}\big).$ 
Using that $\mathcal{F}(H^N_c)= \mathfrak{h}^N_c$ as well as $\mathcal{F}(H^N_r)= \mathfrak{h}^N_r$
it then follows from Theorem \ref{Theorem3.1} that 
\[
z_1(n) = -\hat\varphi_1(-n) + \frac{1}{n^{N+1}}\ell^2_n.
\]
Similarly one shows that $z_2(n)=-\hat\varphi_2(n) + \frac{1}{n^{N+1}}\ell^2_n$.
Going through the arguments of the proof one verifies that these estimates hold locally uniformly on $W\cap H^N_c$ and uniformly on bounded subsets of $W\cap H^N_r.$ This shows that $\Phi-\mathcal{F}$ maps $W\cap H^N_c$ into $\mathfrak{h}^{N+1}_c$ and that the map $\Phi-\mathcal{F}: H^N_r\to \mathfrak{h}^{N+1}_r$ is bounded. As each component of $\Phi-\mathcal{F}$ is analytic on $W$ and $\Phi-\mathcal{F}$ is locally bounded on $W\cap H^N_c$ as a map from $W\cap H^N_c$ to $\mathfrak{h}^{N+1}_c$ it follows that $\Phi-\mathcal{F}: H^N_c\to \mathfrak{h}^{N+1}_c$ is real analytic by \cite{GK}, Theorem A.5. This finishes the proof of Theorem \ref{Theorem1.1}.
\end{proof} 
\begin{proof}
[Proof of Corollary \ref{Corollary1.2}]
 (i) Note that 
$\Phi^{-1}- \mathcal{F}^{-1}= - \mathcal{F}^{-1}\circ \left(\Phi-\mathcal{F}\right)\circ \Phi^{-1}.$
As by \cite{Mo}, $\Phi^{-1}: \mathfrak{h}^N_r \to H^N_r$ is bounded for any $N\in \mathbb{Z}_{\geq 1},$ the stated result follows from Theorem \ref{Theorem1.1}. 
(ii) For any $s\in \mathbb{R}_{\geq 1},$ let $N:=\lfloor s \rfloor$. Then $N\leq s< N+1.$ As by Theorem \ref{Theorem1.1}, the map $\Phi-\mathcal{F}: H^{N}_r\to \mathfrak{h}^{N+1}_r$ is real analytic and bounded and $H_r^s \hookrightarrow H^N_r$ as well as $\mathfrak{h}^{N+1}_r
\hookrightarrow   \mathfrak{h}^{s}_r$ are linear bounded embeddings, $\Phi: H^s_r\to \mathfrak{h}^{s}_r$ is real analytic and bounded. Arguing as in the proof of item (i) it then follows from \[\Phi^{-1}= \mathcal{F}^{-1}+ \left(\Phi^{-1}- \mathcal{F}^{-1}\right)=\mathcal{F}^{-1} - \mathcal{F}^{-1}\circ \left(\Phi-\mathcal{F}\right)\circ \Phi^{-1} \]
that $\Phi^{-1}: \mathfrak{h}^s_r \to H^s_r$  is real analytic and bounded as well. In particular, $\Phi:  H^s_r \to \mathfrak{h}^s_r$ is a diffeomorphism. 
\end{proof}

\begin{proof}
[Proof of Corollary \ref{Corollary1.3}]
For any $s\in \mathbb{R}_{\geq 1},$ let $N:=\lfloor s \rfloor$.
Let $\left(\varphi^{(j)}\right)_{j\geq 1}$
be any sequence in $H^s_r$ which converges weakly to an element  $\varphi\in H^s_r$. As by Rellich's theorem, the embedding $H_r^s\hookrightarrow H_r^0$ is compact, $\varphi^{(j)} \underset{j\to \infty}{\longrightarrow} \varphi$ strongly in $H^0_r$ and hence by Theorem \ref{Theorem2.1}, 
$\Phi\left(\varphi^{(j)}\right) \underset{j\to \infty}{\longrightarrow} \Phi(\varphi)$ strongly in $\mathfrak{h}^0_r.$ In particular, one has $\Phi_n\left(\varphi^{(j)}\right) \underset{j\to \infty}{\longrightarrow} \Phi_n(\varphi)$ for each component $\Phi_n$ of $\Phi= \left(\Phi_n \right)_{n\in \mathbb{Z}}.$
As by Corollary \ref{Corollary1.2} (ii), $\left(\Phi\left(\varphi^{(j)}\right)\right)_{j\geq 0}$  is bounded in $\mathfrak{h}_r^s$ it then follows that $\Phi\left(\varphi^{(j)}\right) \underset{j\to \infty}{\rightharpoondown} \Phi(\varphi)$ weakly in $\mathfrak{h}_r^s.$ Arguing in a similar way one sees that $\Phi^{-1}: \mathfrak{h}^s_r \to H^s_r$ is weakly continuous as well.
\end{proof}



\begin{thebibliography}{123}

\bibitem{Bo0} J. Bourgain: {\em Fourier transform restriction phenomena for certain lattice subsets and applications to 
nonlinear evolution equations: Part I: Schr\"odinger equations}, 
Geom. Funct. Anal., $\bf 3$ (1993), 209-262

\bibitem{GK} B. Gr\'ebert, T. Kappeler: {\em The defocusing NLS equation
and its normal form}, EMS, 2014

\bibitem{KdV&KAM} T. Kappeler, J. P\"oschel: {\em KdV \& KAM}, 
Springer, Berlin, 2003

\bibitem{KST1}  T. Kappeler, B. Schaad, P. Topalov: {\em Qualitative features of periodic solutions of KdV}, 
Comm. Partial Differential Equations, $\bf 38$(2013), no. 9, 1626-1673

\bibitem{KST2}  T. Kappeler, B. Schaad, P. Topalov: {\em Asymptotics of spectral quantities of Zakharov-Shabat operators},
arXiv:1503.04850

\bibitem{KST3}  T. Kappeler, B. Schaad, P. Topalov:  {\em Scattering-like phenomena of the periodic defocusing NLS equation}, preprint
  
\bibitem {Ku} S. Kuksin: {\em Damped-driven KdV and effective equations for long-time behavior of its solution}, 
Geom. Funct. Anal., $\bf 20$(2010), 1431-1463
  
\bibitem {KuPe} S. Kuksin, G. Perelman: {\em Vey theorem in infinite dimensions
and its application to KdV}, Discrete Contin. Dyn. Syst. A, $\bf 27$(2010), no 1, 1-24
  
\bibitem {KuPi} S. Kuksin, A. Piatnitski: {\em Khasminskii-Whitham averaging for randomly perturbated KdV equation}, 
J. Math. Pures Appl., $\bf 89$(2008), no. 4, 400-428
  
\bibitem{Mo} J. Molnar: {\em New estimates of the non-linear Fourier transform for the defocusing NLS equation},  
Int. Math. Res. Notices, $\bf 2014$, doi: 10.1093/im\-rn/rnu208 

\end{thebibliography}
\end{document}